\documentclass[12pt]{amsart}

\usepackage{amsmath,amsfonts,amscd,amsthm,amssymb,latexsym}

\font\teneufm=eufm10 \font\seveneufm=eufm7 \font\fiveeufm=eufm5
\newfam\frakturfam
\textfont\frakturfam=\teneufm \scriptfont\frakturfam=\seveneufm
\scriptscriptfont\frakturfam=\fiveeufm

\newtheorem{pr}{Proposition}

\newtheorem{lm}{Lemma}
\newtheorem{theor}{Theorem}
\newtheorem{co}{Corollary}

\def\bee{\begin{eqnarray}}
\def\bes{\begin{eqnarray*}}
\def\eee{\end{eqnarray}}
\def\ees{\end{eqnarray*}}

\def\Proof{{\sl Proof.}\ }

\title{On the Lie-solvability of Novikov algebras}

\begin{document}
\date{}
\maketitle

\begin{center}

{\bf Kaisar Tulenbaev}\footnote{Department of Mathematics,
 Suleyman Demirel University,
 Almaty, 050040, Kazakhstan,
e-mail: {\em kaisar.tulenbayev@sdu.edu.kz}}, 
{\bf Ualbai Umirbaev}\footnote{Department of Mathematics,
 Wayne State University,
Detroit, MI 48202, USA; Department of Mathematics,
Al-Farabi Kazakh National University, Almaty, 050040, Kazakhstan;
Institute of Mathematics of the SB of RAS, Novosibirsk, 630090, Russia;
and Institute of Mathematics and Mathematical Modeling, Almaty, 050010, Kazakhstan,
e-mail: {\em umirbaev@wayne.edu}}, 
 and 
{\bf Viktor Zhelyabin}\footnote{Institute of Mathematics of the SB of RAS, Novosibirsk, 630090, Russia,
e-mail: {\em vicnic@math.nsc.ru}}

\end{center}

\begin{abstract} 
We prove that any  Novikov algebra over a field of characteristic $\neq 2$ is Lie-solvable if and only if its commutator ideal $[N,N]$ is right nilpotent. We also construct examples of infinite-dimensional Lie-solvable Novikov algebras $N$ with non nilpotent  commutator ideal $[N,N]$. 
\end{abstract}

\noindent
 {\bf Mathematics Subject Classification (2020):} 17D25, 17B30, 17B70

\noindent
{\bf Key words:} Novikov algebra, Lie-solvability, nilpotency

\section{Introduction}

\hspace*{\parindent} 

An algebra  $N$ over a field  $K$ is called a {\em Novikov} algebra if
it satisfies the following  identities:
\bee\label{f1}
(x,y,z)=(y,x,z), 
\eee
\bee\label{f2}
(xy)z=(xz)y,
\eee
 where  $(x,y,z)=(xy)z-x(yz)$ is the associator of elements $x,y,z$.
 
 Recall that any algebra satisfying the identity (\ref{f1}) is called {\em left-symmetric}. Left-symmetric algebras are Lie-admissible, i.e., every left-symmetric algebra $L$ becomes a Lie algebra with respect to the commutator $[x,y]=xy-yx$. This Lie algebra is denoted by $L^{(-)}$ and is called the {\em commutator algebra} of $L$. 
 
 Left-symmetric algebras arise in many areas of mathematics and physics \cite{Burde06}. The defining identities of Novikov 
 algebras first appeared in the study of Hamiltonian operators in the formal calculus of variations by I.M. Gelfand and
I.Ya. Dorfman \cite{GD79}. These identities played a crucial role in the classification of linear Poisson brackets of hydrodynamical type by A.A. Balinskii and  S.P. Novikov \cite{BN85}.  

In 1987 E.I. Zelmanov \cite{Zel} proved that any finite dimensional simple Novikov algebra over a field $K$ of characteristic zero is one-dimensional.  V.T. Filippov \cite{Fil89} constructed a wide class of simple Novikov algebras of characteristic $p\geq 0$. 
 J.M. Osborn \cite{Osborn1992,Osborn1992CommAl,Osborn1994} and X. Xu \cite{Xu96,Xu01} continued the study of simple finite dimensional algebras over fields of positive characteristic and simple infinite dimensional algebras over fields of characteristic zero.  A complete classification of finite dimensional simple Novikov algebras over algebraically closed fields of characteristic $p>2$ is given in \cite{Xu96}. 
   
E.I. Zelmanov also proved that if $N$ is a finite
dimensional right nilpotent Novikov algebra then $N^2$ is nilpotent \cite{Zel}. In 2001 V.T. Filippov \cite{Fil01} proved that any left-nil Novikov algebra of bounded index over a field of characteristic zero is nilpotent. 
A.S. Dzhumadildaev and K.M. Tulenbaev \cite{DzT} proved that any right-nil Novikov algebra of bounded index $n$ is right nilpotent if the characteristic $p$ of the field $K$ is $0$ or $p>n$. In 2020 I. Shestakov and Z. Zhang proved \cite{ShZh} that for any
Novikov algebra $N$ over a field the following conditions are
equivalent:

$(i)$ $N$ is solvable;

$(ii)$ $N^2$ is nilpotent;

$(iii)$ $N$ is right nilpotent.

U.U. Umirbaev and V.N. Zhelyabin proved  \cite{UZh21,ZhU} that any $\mathbb{Z}_n$-graded Novikov algebra with solvable $0$-component is solvable.

It is well known \cite{Jac} that if $L$ is a finite dimensional solvable Lie algebra over a field $K$ of characteristic zero then $[L,L]$ is nilpotent. 
In 1973 Yu.P. Razmyslov proved \cite{Raz73} that over a field $K$ of characteristic zero $[L,L]$ is nilpotent for any algebra $L$  from any proper subvariety of the variety of algebras generated by the simple three dimensional Lie algebra $\mathrm{sl}_2(K)$.  There was a long standing conjecture about solvable algebras of the variety of algebras generated by the Witt algebra $W_1$. 

{\em Conjecture 1}. If $L$ is a solvable algebra of the variety of algebras generated by the Witt algebra $W_1$, then is it true 
that $[L,L]$ is nilpotent? 

This conjecture was proven to be not true by A. Mishchenko \cite{Mish} in 1988.

The variety of Lie algebras generated by the Witt algebra $W_1$ is closely related to the variety of Novikov algebras. 
Let $K[x]$ be the algebra of all polynomials in one variable $x$ over a field $K$. Consider $K[x]$ as a differential algebra with derivation $\partial=\frac{\partial}{\partial x}$. Then $K[x]$ is a simple differential algebra over a field of characteristic zero. 
With respect to the product 
\bes
f\circ g= fg'
\ees
the vector space $K[x]$ becomes a Novikov algebra. We denote this algebra by $L_1$. The construction described above is called the {\em Gelfand-Dorfman} construction for Novikov algebras.  Recently, L.A. Bokut, Y. Chen, and Z. Zhang \cite{BCZ18} proved that any Novikov algebra over a field of characteristic zero is a subalgebra of a Novikov algebra obtained from some differential algebra by the Gelfand-Dorfman construction. 

Notice that $K[x]$ becomes a Lie algebra with respect to the product 
\bes
[f,g]=fg'-gf'. 
\ees
This algebra is a well known Witt algebra $W_1$. This construction of Lie algebras is also studied by many specialists \cite{PR16,Poin18,Raz86}.
In this case the differential enveloping algebra of a Lie algebra is called the {\em Wronskian} enveloping algebra \cite{Poin18}. Although there are many interesting results, the class of Lie algebras embeddable into their Wronskian enveloping algebras is not described yet. 

{\em Conjecture 2}. A Lie algebra over a field of characteristic zero is embeddable into its Wronskian enveloping algebra if and only if it belongs to the variety of algebras generated by the Witt algebra $W_1$. 

Notice that the commutator algebra of $L_1$ is the Witt algebra $W_1$. For this reason we call $L_1$ the {\em Novikov-Witt} algebra \cite{KUZ21}. This is the first algebra in the list of left-symmetric Witt algebras $L_n$ \cite{U16}. 
The variety of Novikov algebras is generated by the Novikov-Witt algebra $L_1$  in characteristic zero \cite{MLU11N}. The identities of the Witt algebras $W_n$ are studied mainly by Yu.P. Razmyslov \cite{Raz94} and the identities of the left-symmetric Witt algebras $L_n$ are studied in \cite{KU16}

We say that a Novikov algebra $N$ is {\em Lie-solvable} if the Lie algebra $N^{(-)}$ is solvable.  
It is known that every finite dimensional Novikov algebra over a field is Lie-solvable \cite{BD06}. 
Recently 
Z. Zhang and T.G. Nam \cite{ZhNam} proved that if a Novikov algebra is Lie-nilpotent then 
its ideal generated by all commutators $[a,b]$ is nilpotent.

This paper is devoted to the study of Lie-solvable Novikov algebras. We noticed that the space of commutators $[N,N]$ of a Novikov algebra $N$ is an ideal of $N$ over a field of characteristic $\neq 2$. We prove that a Novikov algebra $N$ over a field of characteristic $\neq 2$ is Lie-solvable if and only if $[N,N]$ is right nilpotent. Using Mishchenko's example \cite{Mish}, we constructed examples of Lie-solvable Novikov algebras with non nilpotent $[N,N]$. 

The right nilpotency of $[N,N]$ for Lie-solvable Novikov algebras means that Conjecture 1 was not baseless. This property just cannot be expressed in the language of Lie algebras. Notice that if $[N,N]$ is right nilpotent then  $[N,N] [N,N]$ is nilpotent by the above mentioned result of  I. Shestakov and Z. Zhang \cite{ShZh} . This fact suggests us to formulate the following weaker version of Conjecture 1. 

{\em Conjecture 3}. If $L$ is a solvable algebra of the variety of algebras generated by the Witt algebra $W_1$, then is it true 
that $[[L,L],[L,L]]$ is nilpotent?

The paper is organized as follows. In Section 2 we give some identities, construction of ideals,  and recall some definitions.  Sections 3 is devoted to the proof of the main result on the right nilpotency of $[N,N]$. Examples of Novikov algebras $N$ with non nilpotent $[N,N]$ are given in Section 4.

\section{Identities, ideals, and some definitions}

\hspace*{\parindent}

As we mentioned above, any left -symmetric algebra is Lie-admissible, i.e., satisfies the Jacobi identity 
\bee\label{f3}
[[x,y],z]+[[y,z],x]+[[z,x],y]=0. 
\eee
 Moreover, in the class of Novikov algebras this identity splits into the identities 
\bee\label{f4}
[x,y]z+[y,z]x+[z,x]y=0
\eee
and 
\bee\label{f5}
x[y,z]+y[z,x]+z[x,y]=0.
\eee
Indeed, 
\bes
[x,y]z+[y,z]x+[z,x]y=(xy)z-(yx)z+(yz)x-(zy)x+(zx)y-(xz)y=0
\ees
by (\ref{f2}). This proves  (\ref{f4}). Using (\ref{f3}) and (\ref{f4}) we also get  (\ref{f5}). 

It is useful to write the identity (\ref{f2}) in the form
\bee\label{f6}
x[y,z]=(x,z,y)-(x,y,z). 
\eee

Any nonassociative algebra satisfies (see \cite{ZSSS}) the identity
\bes
[xy,z]-x[y,z]-[x,z]y=(x,y,z)-(x,z,y)+(z,x,y).
\ees
Using this and (\ref{f6}) we get 
\bee\label{f7}
(z,x,y)=[xy,z]-[x,z]y.
\eee

 The identities (\ref{f1}) and (\ref{f2}) easily imply that
\bee\label{f8}
(xy,z,t)=(x,z,t)y
\eee
and, consequently, 
\bee\label{f9}
(x,yz,t)=(x,y,t)z.
\eee

Recall that any nonassociative algebra also satisfies (see \cite{ZSSS}) the identity
\bee\label{f10}
(x,y,zt)=(x,yz,t)-(xy,z,t)+x(y,z,t)+(x,y,z)t,
\eee

Then (\ref{f8}), (\ref{f9}), and (\ref{f10}) give that
\bee\label{f11}
(x,y,zt)=(x,y,t)z-(x,z,t)y+x(y,z,t)+(x,y,z)t.
\eee

It is well known that if $I$ and $J$ are ideals of a Novikov algebra $N$, then $IJ$ is  an ideal of $N$.
\begin{lm}\label{l1} \cite{ZhT08} 
In any Novikov algebra $N$ the space of associators $(N,N,N)$ is an ideal of $N$. 
\end{lm}
\Proof The space $(N,N,N)$  is a right ideal by (\ref{f8}) or (\ref{f9}).  Applying (\ref{f10}) or (\ref{f11}), we get that  $(N,N,N)$ is also a left ideal. 
$\Box$

\begin{lm}\label{l2} Any Novikov algebra over a field of characteristic $\neq 2$ satisfies the identities

\bee\label{f12}
(a,b,x)=\frac{1}{2}([ax,b]-[a,bx]), 
\eee
\bee\label{f13}
[a,b]x=\frac{1}{2}([ax,b]+[a,bx]),
\eee
and 
\bee\label{f14}
x[a,b]=[[x,a],b]+[a,[x,b]]+\frac{1}{2}([ax,b]+[a,bx]). 
\eee\end{lm}
\Proof Applying once the identity (\ref{f2}), we get 
\bes
(a,b,x)=(ab)x-a(bx)=(ab)x-[a,bx]-(bx)a\\
=(ab)x-[a,bx]-(ba)x=[a,b]x-[a,bx]. 
\ees
Consequently, 
\bes
(b,a,x)=[b,a]x-[b,ax]. 
\ees
By (\ref{f1}), we get 
\bes
2(a,b,x)=[a,b]x-[a,bx]+[b,a]x-[b,ax]
\ees
and 
\bes
[a,b]x-[a,bx]=[b,a]x-[b,ax].
\ees
Consequently,
\bes
2(a,b,x)=[ax,b]-[a,bx]
\ees
and 
\bes
2[a,b]x=[ax,b]+[a,bx], 
\ees
which imply (\ref{f12}) and (\ref{f13}), respectively

Using (\ref{f3})  and (\ref{f13}), we get 
\bes
x[a,b]=[x,[a,b]]+[a,b]x=[[x,a],b]+[a,[x,b]]+\frac{1}{2}([ax,b]+[a,bx]),  
\ees
i.e., (\ref{f14}) holds. $\Box$

\begin{co}\label{c1} Let  $N$ be a Novikov algebra over a field of characteristic $\neq 2$.  Then the following statements are true: 

$(i)$ If $I$ and $J$ are right ideals of $N$ then $[I,J]$ is a right ideal of $N$; 

 $(i)$ If $I$ and $J$ are ideals of $N$ then $[I,J]$ is an ideal of $N$.
\end{co}

At the end of this section we recall the definitions of solvable, nilpotent, and right nilpotent algebras.

Let $A$ be an arbitrary algebra. The powers of $A$ are defined inductively by $A^1=A$ and
\bes
A^m=\sum_{i=1}^{m-1}A^{i}A^{m-i}
\ees
 for all positive integers $m\geq 2$. The algebra $A$ is called {\em nilpotent}
if  $A^{m}=0$ for some positive integer  $m$.

The right powers of $A$ are defined inductively by $A^{[1]}=A$ and $A^{[m+1]}=A^{[m]}A$ for all integers $m\geq 1$. The algebra $A$ is called {\em right nilpotent} if there exists a positive integer $m$ such that $A^{[m]}=0$. In general, the right nilpotency of an algebra does not imply its nilpotency. This is also true in the case of Novikov algebras.

{\em Example 1}. \cite{Zel} Let  $N=Fa+Fb$ be a vector space of dimension 2. The product on $N$ is defined as
 $$ab=b,a^2=b^2=ba=0.$$
 It is easy to check that $N$ is a right nilpotent Novikov algebra, but not nilpotent.

The derived powers of $A$ are defined by $A^{(0)}=A$,
$A^{(1)}=A^2$, and $A^{(m)}=A^{(m-1)}A^{(m-1)}$ for all positive
integers $m\geq 2$. The algebra $A$ is called {\em solvable} if
$A^{(m)}=0$ for some positive integer  $m$. Every right nilpotent
algebra is solvable, and, in general, the converse is not true. But
every solvable Novikov algebra is right
nilpotent \cite{ShZh}.

A Novikov algebra $N$ is called {\em Lie-solvable} if the Lie algebra $N^{(-)}$ is solvable.

\section{Lie-solvable Novikov algebras}

\hspace*{\parindent}

A Novikov algebra $N$ is called {\em Lie-metabelian} if it satisfies the identity
\bee\label{g1}
[[x,y],[z,t]]=0.
\eee

In any algebra we denote by $x_1x_2\ldots x_k$ the right normed product $(\ldots(x_1x_2)\ldots )x_k$ of elements
 $x_1,x_2,\ldots,x_k$. 

\begin{lm}\label{l3} Any Lie-metabelian Novikov algebra $N$ over a field $K$ of characteristic $\neq 2$ satisfies the identity  
\bee\label{g2} 
([x,y],[z,t],s)=0.
\eee
\end{lm}
\Proof Using (\ref{g1}) and (\ref{f13}), we immediately get 
\bes
[[a,b]x,[y,z]]=0.
\ees
Then the identities (\ref{g1}) and  (\ref{f2}) imply that
\bes
([y,z],[a,b],x)=[y,z][a,b]x-[y,z]([a,b]x)=[a,b][y,z]x-[y,z]([a,b]x)\\
=[a,b]x[y,z]-[y,z]([a,b]x)=[[a,b]x,[y,z]]=0. \ \ \Box
\ees

\begin{co}\label{c2} 
The ideals $(N,N, N)$ and $[N,N]$  of a Novikov algebra $N$ over a field $K$ of characteristic $\neq 2$ are associative and commutative and $(N,N, N)\subseteq [N,N]$.
\end{co}
\Proof Notice that $(N,N, N)$ is an ideal of $N$ by Lemma \ref{l1} and $[N,N]$ is an ideal of $N$ by Corollary \ref{c1}. The identity (\ref{f12}) implies that $(N,N, N)\subseteq [N,N]$. The identities (\ref{g1}) and (\ref{g2}) imply that $[N,N]$ is an associative and commutative algebra. $\Box$

\begin{lm}\label{l4} Any Lie-metabelian Novikov algebra $N$ over a field $K$ of characteristic $\neq 2$ satisfies the identities 
\bee\label{g3} 
(x,[y,z],t)[a,b]=0, 
\eee
\bee\label{g4} 
([x,y][z,t],a,b)=0, 
\eee
and
\bee\label{g5}
[x,y][z,t](a,b,c)=0. 
\eee
\end{lm}
\Proof The identities  (\ref{f14})  and (\ref{g2}) imply that 
\bes
(x[a,b],[y,z],t)=0, 
\ees
Then, by (\ref{f8}), 
\bes
(x,[y,z],t)[a,b]=(x[a,b],[y,z],t)=0, 
\ees 
i.e., (\ref{g3}) holds. 
Using  (\ref{f8}),  (\ref{f1}), and  (\ref{g3}), we get 
\bes
([x,y][z,t],a,b)=([x,y],a,b)[z,t]=(a,[x,y],b)[z,t]=0,
\ees
i.e., (\ref{g4}) also holds. 
By (\ref{f10}), 
\bes
[x,y][z,t](a,b,c)=([x,y][z,t],a,bc)\\
-([x,y][z,t],ab,c)+([x,y][z,t]a,b,c)-([x,y][z,t],a,b)c.
\ees
Using (\ref{g4}), from this we get (\ref{g5}). 
$\Box$ 

\begin{lm}\label{l5} Let $N$ be a Lie-metabelian Novikov algebra over a field $K$ of characteristic $\neq 2$. Then $(N,N,N)^3=[N,N]^4=0$. 
\end{lm}
\Proof The identity (\ref{g5}) implies that $(N,N,N)^2(N,N,N)=0$ since $(N,N,N)\subseteq [N,N]$ by Corollary \ref{c2}. Consequently, $(N,N,N)^3=0$ since $(N,N,N)$ is 
associative. 

Notice that $N[N,N]\subseteq (N,N,N)$ by (\ref{f6}). Consequently, $[N,N]^2\subseteq (N,N,N)$. Then (\ref{g5}) implies that
$[N,N]^2[N,N]^2=0$. This gives $[N,N]^4=0$ since $[N,N]$ is an associative algebra. $\Box$

\begin{theor}\label{t1}
Let $N$ be a Lie-solvable Novikov algebra over a field of characteristic $\neq 2$. Then the ideal
$[N,N]$ is right nilpotent.
\end{theor} 
\Proof 
Let $N$ be a Lie-solvable Novikov algebra with Lie-solvable index $n$. We prove the statement of the theorem by induction on $n$. By Lemma \ref{l5}, this is true for $n=2$. Suppose that $n\geq 3$. Then $[N,N]$ is a Lie-solvable Novikov algebra with Lie-solvable index $n-1$. 

By the induction hypothesis $[[N,N],[N,N]]$ is a right nilpotent ideal of $N$. Notice that $[N,N]^4\subseteq [[N,N],[N,N]]$. 
Consequently, $[N,N]$ is a solvable ideal of $N$. Recall that every solvable Novikov algebra is right nilpotent \cite{ShZh}. 
Therefore $[N,N]$ is a right nilpotent ideal of $N$. $\Box$

\section{Lie-solvable Novikov algebras with non-nilpotent commutator ideal}

\hspace*{\parindent}

Let $K[x]$ be the polynomial algebra over a field $K$ of characteristic zero in one variable $x$. Recall that the Witt algebra $W_1$ is the Lie algebra of all derivations of $K[x]$. Any element of  $W_1$ can be written in the form 
\bes
f\partial, 
\ees
where $f\in K[x]$ and $\partial=\frac{\partial}{\partial x}$. 
The vector space of $W_1$ with respect to the product 
\bes
f\partial \circ g\partial=fg'\partial
\ees
becomes a Novikov algebra \cite{GD79}. This algebra is denoted by $L_1$ and is called the {\em Novikov-Witt} algebra \cite{KUZ21}. The elements 
\bes
e_n=x^{n+1}\partial, n\geq -1, 
\ees
form a linear basis of $L_1$ and 
\bes
e_i\circ e_j=(j+1) e_{i+j}
\ees
for all $i,j\geq -1$. Consequently, 
\bes
L_1=\oplus_{i\geq -1} Ke_i 
\ees
is a graded algebra. 

Set 
\bes
R=Ke_{-1}\oplus Ke_0. 
\ees
Notice that $R$ is a subalgebra of $L_1$. The left and right actions of elements of $R$ on  $L_1$ are naturally defined 
since $R$ is a subalgebra of $L_1$. We denote an isomorphic copy of this $R$-bimodule $L_1$ by $M$ and assume $e_i$ corresponds to $f_i\in M$ for all $i$.  This means that 
\bes
M=\oplus_{i\geq -1} Kf_i 
\ees
and 
\bes
e_i\circ f_j=(j+1) f_{i+j}, \ \ f_j\circ e_i=(i+1) f_{j+i}
\ees
for all $i=-1,0$ and $j\geq -1$.

Since $R$ is a subalgebra of $L_1$ it follows that $L_1$ is a {\em Novikov bimodule} over $R$, i.e., $M$ is a Novikov $R$-bimodule. By definition this means that the space 
\bes
N=R\oplus M 
\ees
with the product 
\bes
(r_1+m_1)(r_2+m_2)=r_1\circ r_2+r_1\circ m_2+m_1\circ r_2,    
\ees
for all $r_1,r_2\in R$ and $m_1,m_2\in M$, is a Novikov algebra. Recall that $N$ is called the {\em zero split extension} of $R$ by $M$. 
\begin{pr}\label{p1}
The Novikov algebra $N$ is Lie-solvable of index $3$ over a field of characteristic zero and $[N,N]$ is not nilpotent. 
\end{pr}
\Proof Obviously, 
\bes
[N,N]=[R,R]\oplus [R,M]
\ees
and $[R,R]=Ke_{-1}$. Moreover,  
\bes
[R,M]=M 
\ees
since 
\bes
[e_{-1},f_j]=(j+1)f_{j-1}  
\ees
 for all $j$. Consequently,  
\bes
[N,N]=Ke_{-1}\oplus M. 
\ees
Obviously, $[N,N]$ is not left nilpotent since $e_{-1}M=M$.  
Furthermore, 
\bes
[[N,N],[N,N]]= M
\ees
and, consequently, $N$ is Lie-solvable of index $3$. $\Box$

The Lie algebra $N^{(-)}$ coincides with  Mishchenko's example from \cite{Mish}. 

Notice that $[N,N]$ is not nilpotent over fields of positive characteristic. In order to adopt this example to the case of positive characteristic, we consider another basis 
\bes
E_i=\frac{1}{(i+1)!}x^{i+1}\partial
\ees
of the space of $L_1$. Recall that binomial coefficients are defined by 
\bes
\binom{n}{k} = \frac{n!}{k!(n-k)!}
\ees
for all integers $n\geq k\geq 0$. For conveniency of notation we set $\binom{n}{k} =0$ if $n<k$. Then  
\bee\label{h1}
E_i\circ E_j=\binom{i+j+1}{i+1} E_{i+j}
\eee
for all $i,j\geq -1$. 

Denote by $L$ the abstract algebra over a field $K$ of arbitrary characteristic with a linear basis 
\bee\label{h2}
E_{-1},E_0,E_1,\ldots,E_k,\ldots
\eee
 and with multiplication defined by (\ref{h1}). 
\begin{lm}\label{}
The algebra $L$ is a Novikov algebra. 
\end{lm}
\Proof Let $S$ be the free $\mathbb{Z}$-module with a free basis (\ref{h2}). We turn $S$ into a $\mathbb{Z}$-algebra by (\ref{h1}). If the characteristic of $K$ is zero, then 
\bes
S_K=S\otimes_{\mathbb{Z}}K
\ees
 is a free $K$-module with a linear basis (\ref{h2}) since $S$ and $K$ are both free $Z$-modules. Consequently, $S$ imbeds into $S\otimes_{\mathbb{Z}}K$. Then $K$-algebras $L$, $S_K$, and $L_1$ are isomorphic by construction. Consequently, 
$L$ is a Novikov algebra over $K$ and  $S$ is a Novikov  algebra over $\mathbb{Z}$. 

Assume that the characteristic of the field $K$ is $p>0$.  Obviously, 
\bes
S_1=S\otimes_{\mathbb{Z}}(\mathbb{Z}/p\mathbb{Z})=S/(p\mathbb{Z})S
\ees
is a free $\mathbb{Z}/p\mathbb{Z}$-module with a basis (\ref{h2}). 
Since $S$ is a Novikov  algebra over $\mathbb{Z}$ it follows that $S_1$ is a Novikov algebra over $\mathbb{Z}/p\mathbb{Z}$. 
This implies that 
\bes
S_K=S\otimes_{\mathbb{Z}}K=S_1\otimes_{\mathbb{Z}/p\mathbb{Z}}K
\ees
is a Novikov algebra over $K$ with a linear basis  (\ref{h2}). Obviously, $L\simeq S_K$. $\Box$

Let 
\bes
R=KE_{-1}\oplus KE_0. 
\ees
be the two dimensional subalgebra of $L$. Consider $L$ as an $R$-bimodule. Denote by $M$ an isomorphic copy of $R$-bimodule $L$  and denote by $F_i$ the images of $E_i$ in $M$ for all $i$.  Then 
\bes
M=\oplus_{i\geq -1} KF_i 
\ees
and 
\bes
E_i\circ F_j=\binom{i+j+1}{i+1} F_{i+j}, \ \ F_j\circ E_i=\binom{i+j+1}{j+1} F_{j+i}
\ees
for all $i=-1,0$ and $j\geq -1$.

Then the zero split extension 
\bes
N=R\oplus M 
\ees
of $R$ by $M$ is a Novikov algebra. 

\begin{pr}\label{p2}
The Novikov algebra $N$ is Lie-solvable of index $3$ over an arbitrary field $K$ and $[N,N]$ is not nilpotent. 
\end{pr}
\Proof Obviously, $[R,R]=KE_{-1}$ and $[R,M]=M$ since $E_{-1}\circ F_j=F_{j-1}$ for all $j\geq -1$. Consequently, 
\bes
[N,N]=KE_{-1}\oplus M. 
\ees
Obviously, $[N,N]$ is not left nilpotent since $E_{-1}M=M$.  
Furthermore, 
\bes
[[N,N],[N,N]]= M
\ees
and, consequently, $N$ is Lie-solvable of index $3$. $\Box$

\begin{center}
{\large Acknowledgments}
\end{center}

\hspace*{\parindent}

This research is supported by the Russian Science Foundation (project
21-11-00286) and by the grants of the Ministry of Education and Science
of the Republic of Kazakhstan (projects AP08855944 and  AP09261086).


\begin{thebibliography}{99}



\bibitem{BN85}  I.M.  Balinskii, S.P. Novikov,   Poisson brackets of hydrodynamic type, Frobenius algebras and Lie
algebras. (Russian) Dokl. Akad. Nauk SSSR 283 (1985), no. 5, 1036--1039.


\bibitem{BCZ18} L.A. Bokut, Y. Chen, Z. Zhang, On free Gelfand-Dorfman-Novikov-Poisson algebras and a PBW theorem. J. Algebra 500 (2018), 153--170.

\bibitem{Burde06} D. Burde, Left-symmetric algebras, or pre-Lie algebras in geometry and physics, Cent. Eur. J. Math. 4 (2006), no. 3, 323--357.

 \bibitem{BD06}  D. Burde, K. Dekimpe, Novikov structures on solvable Lie algebras. 
J. Geom. Phys. 56 (2006), no. 9, 1837--1855. 



\bibitem{DzT} A.S. Dzhumadil'daev, K.M. Tulenbaev, Engel theorem for Novikov algebras. Comm. Algebra 34 (2006),  no. 3, 883--888.


\bibitem{Fil89} V.T. Filippov,
A class of simple nonassociative algebras. (Russian)
Mat. Zametki 45 (1989), no. 1, 101--105; translation in
Math. Notes 45 (1989), no. 1--2, 68--.


\bibitem{Fil01} V.T. Filippov, On right-symmetric and Novikov nil algebras of bounded index. (Russian) Mat. Zametki 70
 (2001), no. 2, 289--295; translation in Math. Notes 70 (2001), no. 1--2,
 258--263.




\bibitem{GD79} I.M. Gel'fand,  I.Ya. Dorfman,   Hamiltonian operators and algebraic structures related to
them. (Russian) Funktsional. Anal. i Prilozhen 13  (1979), no. 4, 3--30.



\bibitem{Jac} N. Jacobson, Lie algebras. Republication of the 1962 original. Dover Publications, Inc., New York, 1979. ix+331 pp. 



\bibitem{KUZ21} D. Kozybaev, U. Umirbaev, V. Zhelyabin, Some examples of nonassociative coalgebras and supercoalgebras. 
Linear Algebra Appl. (submitted)

\bibitem{KU16} D. Kozybaev, U. Umirbaev, Identities of the left-symmetric Witt algebras. Internat. J. Algebra Comput. 26 (2016), no. 2, 435--450.




\bibitem{MLU11N} L. Makar-Limanov, U. Umirbaev, The Freiheitssatz for Novikov algebras. TWMS J. Pure Appl. Math. 2 (2011), no. 2, 228--235.

\bibitem{Mish} S.P. Mishchenko, Solvable subvarieties of a variety generated by a Witt algebra. (Russian) Mat. Sb. (N.S.) 136(178) (1988), no. 3, 413--425, 431--432; translation in Math. USSR-Sb. 64 (1989), no. 2, 415--426. 

\bibitem{Osborn1992} J.M. Osborn, Novikov algebras. Nova J. Algebra Geom. 1 (1992), no. 1, 1--13.

\bibitem{Osborn1992CommAl} J.M. Osborn, Simple Novikov algebras with an idempotent. Comm. Algebra 20 (1992), no. 9, 2729--2753.

\bibitem{Osborn1994} J.M. Osborn, Infinite-dimensional Novikov algebras of characteristic $0$. J. Algebra 167 (1994), no. 1, 146--167.

\bibitem{PR16} G. Pogudin, Yu.P. Razmyslov, Prime Lie algebras satisfying the standard Lie identity of degree 5. J. Algebra 468 (2016), 182--192.


\bibitem{Poin18} L. Poinsot, The solution to the embedding problem of a (differential) Lie algebra into its Wronskian envelope. Comm. Algebra 46 (2018), no. 4, 1641--1667. 

\bibitem{Raz73} Yu. P. Razmyslov, 
The existence of a finite basis for the identities of the matrix algebra of order two over a field of characteristic zero. (Russian)
Algebra i Logika 12 (1973), 83--113.


\bibitem{Raz86} Yu.P. Razmyslov, Simple Lie algebras satisfying the standard Lie identity of degree 5. (Russian)
Izv. Akad. Nauk SSSR Ser. Mat. 49 (1985), no. 3, 592--634.

\bibitem{Raz94} Yu.P. Razmyslov, Identities of algebras and their representations. Translated from the 1989 Russian original by A. M. Shtern. Translations of Mathematical Monographs, 138. American Mathematical Society, Providence, RI, 1994. xiv+318 pp.

\bibitem{ShZh} I. Shestakov and Z. Zhang, Solvability and nilpotency of Novikov algebras. Comm. Algebra 48 (2020),  no. 12, 5412--5420.

\bibitem{U16} U.U. Umirbaev, Associative, Lie, and left-symmetric algebras of derivations. Transform. Groups 21 (2016), no. 3, 851--869. 

\bibitem{UZh21} U. Umirbaev, V. Zhelyabin, On the solvability of graded Novikov algebras. Internat. J. Algebra Comput. 31 (2021), no. 7, 1405--1418.

\bibitem{Xu96} X. Xu,
On simple Novikov algebras and their irreducible modules.
J. Algebra 185 (1996), no. 3, 905--934.

\bibitem{Xu01} X. Xu, Classification of simple Novikov algebras and their irreducible modules of characteristic 0.
J. Algebra 246 (2001), no. 2, 673--707.


\bibitem{Zel} E.I. Zel'manov, A class of local translation-invariant Lie algebras.  Dokl. Akad. Nauk SSSR 292 (1987),  no. 6, 1294--1297.

\bibitem{ZhT08} V.N. Zhelyabin, Tikhov, A. S., Novikov-Poisson algebras and associative commutative derivation algebras. (Russian. Russian summary)
Algebra Logika 47 (2008), no. 2, 186--202; translation in
Algebra Logic 47 (2008), no. 2, 107--117. 

\bibitem{ZhU} V. Zhelyabin, U. Umirbaev, On the Solvability of $\mathbb{Z}_3$-Graded Novikov
Algebras. Symmetry  312(2) (2021), 13.

\bibitem{ZSSS} K.A. Zhevlakov, A. M. Slinko, I. P. Shestakov, A. I. Shirshov, Rings
that are Nearly Associative. Academic Press, New York, 1982.

\bibitem{ZhNam} Z. Zhang, T.G. Nam. Lie nilpotent Novikov algebras and Lie solvable Leavitt path algebras.

\end{thebibliography}
\end{document}